# Risk bounds for the non-parametric estimation of Lévy processes

## José E. Figueroa-López[1] and Christian Houdré[2]

*Purdue University and Georgia Institute of Technology*

**Abstract:** Estimation methods for the Lévy density of a Lévy process are developed under mild qualitative assumptions. A classical model selection approach made up of two steps is studied. The first step consists in the selection of a good estimator, from an approximating (finite-dimensional) linear model $\mathcal{S}$ for the true Lévy density. The second is a data-driven selection of a linear model $\mathcal{S}$, among a given collection $\{\mathcal{S}_m\}_{m \in \mathcal{M}}$, that approximately realizes the best trade-off between the error of estimation within $\mathcal{S}$ and the error incurred when approximating the true Lévy density by the linear model $\mathcal{S}$. Using recent concentration inequalities for functionals of Poisson integrals, a bound for the risk of estimation is obtained. As a byproduct, oracle inequalities and long-run asymptotics for spline estimators are derived. Even though the resulting underlying statistics are based on continuous time observations of the process, approximations based on high-frequency discrete-data can be easily devised.

## 1. Introduction

Lévy processes are central to the classical theory of stochastic processes, not only as discontinuous generalizations of Brownian motion, but also as prototypical Markov processes and semimartingales (see [27] and [5] for monographs on these topics). In recent years, continuous-time models driven by Lévy processes have received a great deal of attention mainly because of their applications in the area of mathematical finance (see e.g. [14] and references therein). The scope of these models goes from simple exponential Lévy models (e.g. [2, 10, 12] and [16]), where the underlying source of randomness in the Black-Scholes model is replaced by a Lévy process, to exponential time-changed Lévy processes (e.g. [11]-[13]) and to stochastic differential equations driven by multivariate Lévy processes (e.g. [3, 29]). Exponential Lévy models have proved successful to account for several empirical features observed in time series of financial returns such as heavy tails, high-kurtosis, and asymmetry (see, for example, [10] and [16]). Lévy processes, as models capturing the most basic features of returns and as "first-order approximations" to other more accurate models, should be considered first in developing and testing a successful statistical methodology. However, even in such parsimonious models, there are several issues in performing statistical inference by standard likelihood-based methods.

Lévy processes are determined by three "parameters": a non-negative real $\sigma^2$, a real $\mu$, and a measure $\nu$ on $\mathbb{R} \backslash \{0\}$. These three parameters characterize a Lévy process $\{X(t)\}_{t \geq 0}$ as the superposition of a Brownian motion with drift, $\sigma B(t) + \mu t$, and an independent pure-jump Lévy process, whose jump behavior is specified by the measure $\nu$ in that for any $A \in \mathcal{B}(\mathbb{R})$, whose indicator $\chi_A$ vanishes in a neighborhood of the origin,

$$\nu(A) = \frac{1}{t}\mathbb{E}\left[\sum_{s \leq t}\chi_A\left(\Delta X(s)\right)\right],$$

for any $t > 0$ (see Section 19 of [27]). Here, $\Delta X(t) \equiv X(t) - X(t^-)$ denotes the jump of $X$ at time $t$. Thus, $\nu(A)$ gives the average number of jumps (per unit time) whose magnitudes fall in the set $A$. A common assumption in Lévy-based financial models is that $\nu$ is determined by a function $p : \mathbb{R}\backslash\{0\} \to [0,\infty)$, called the *Lévy density*, as follows

$$\nu(A) = \int_A p(x)dx, \;\; \forall A \in \mathcal{B}(\mathbb{R}\backslash\{0\}).$$

[1]Department of Mathematics, Purdue University, West Lafayette, IN 47907, USA, e-mail: jfiguero@math.purdue.edu

[2]School of Mathematics, Georgia Institute of Technology, Atlanta, GA 30332-0160, USA, e-mail: houdre@math.gatech.edu

Intuitively, the value of $p$ at $x_0$ provides information on the frequency of jumps with sizes "close" to $x_0$.

Estimating the Lévy density poses a nontrivial problem, even when $p$ takes simple parametric forms. Parsimonious Lévy densities usually produce not only intractable marginal densities, but sometimes marginal densities which are not even expressible in a closed form. The current practice of estimation relies on numerical approximations of the density function of $X(t)$ using *inversion formulas* combined with maximum likelihood estimation (see for instance [10]). Such approximations make the estimation computationally expensive and particularly susceptible to numerical errors and mis-specifications. Even in the case of closed form marginal densities, maximum-likelihood based methods present serious numerical problems. For instance, analyzing generalized hyperbolic Lévy processes, the author of [24] notices that the likelihood function is highly flat for a wide range of parameters and good starting values as well as convergence are critical. Also, the separation of parameters and identification between different subclasses is difficult. These issues worsen when dealing with "high-frequency" data. Other calibration methods include methods based on moments, simulation based methods, and multinomial log likelihoods (see e.g. [29] and [6] and references therein). However, our goal in the present paper is not to match the precision of some of these parametric methods, but rather gain in robustness and efficiency using *non-parametric methods*. That is to say, assuming only qualitative information on the Lévy density, we develop estimation schemes for the Lévy density $p$ that provide fairly general function estimators $\hat{p}$.

We follow the so-called *model selection methodology* developed in the context of density estimation in [8], and recently extended to the estimation of intensity functions for Poisson processes in [25]. The essence of this approach is to approximate an infinite-dimensional, nonparametric model by a sequence of finite-dimensional models. This strategy has its origins in *Grenander's method of sieves* (see [17]). Concretely, the procedure addresses two problems. First, the selection of a good estimator $\hat{p}_{_\mathcal{S}}$, called the *projection estimator*, out of an approximating (finite-dimensional) linear model $\mathcal{S}$ for the Lévy density. Second, the selection of a linear model $\mathcal{S}_{\hat{m}}$, among a given collection of linear models $\{\mathcal{S}_m\}_m$, that approximately realizes the best trade-off between the error of estimation from the first step, and the error incurred when approximating the unknown Lévy density by the linear model $\mathcal{S}$. The technique used in the second step has the general flavor of cross-validation via a penalization term, leading to *penalized projection estimators* $\tilde{p}$ (p.p.e.).

Comparing our approach to other non-parametric methods for non-homogeneous Poisson processes (see e.g. [20, 21] and [25]), we will see that the main difficulty here is the fact that the jump process associated with a Lévy process has potentially infinitely many small jumps. To overcome this problem, we introduce a reference measure and estimate instead the Lévy density with respect to this measure. In contrast to [25], our treatment does not rely on the finiteness of such a reference measure. Our main objective here is to estimate the order of magnitude of the mean-square error, $\mathbb{E}\|p - \tilde{p}\|^2$, between the true Lévy density and the p.p.e. To accomplish this, we apply concentration inequalities for functionals of Poisson point processes such as functions of stochastic Poisson integrals (see e.g. [9, 18, 25]). This important statistical application of concentration inequalities is well-known in other contexts such as regression and density estimation (see [8] and references therein). The bound for the risk of estimation leads in turn to *oracle inequalities* implying that the p.p.e. is at least as good (in terms of the long term rate of convergence) as the best projection estimator (see Section 4 for details). Also, combining the bound with results on the approximation of smooth functions by *sieves*, one can determine the long-term rate of convergence of the p.p.e. on certain well-known approximating spaces of functions such as splines.

The statistics underlying our estimators are expressed in terms of deterministic functions of the jumps of the process, and thus, they are intrinsically based on a continuous-time observation of the process during some time period $[0, T]$. Even though this observation scheme has an obvious drawback, statistical analysis under it presents a lot of interest for two reasons. First, very powerful theoretical results can be obtained, thus providing benchmarks of what can be achieved by discrete-data-based statistical methods. Second, since the path of the process can in principle be approximated by high-frequency sampling, it is possible to construct feasible estimators by approximating the continuous-time based statistics using discrete-observations. We use this last idea to obtain estimators by replacing the jumps by increments, based on equally spaced observations of the process.

Let us describe the outline of the paper. We develop the model selection approach in Sections 2 and 3. We proceed to obtain in Section 4 bounds for the risk of estimation, and consequently prove *oracle inequalities*. In Section 5 the rate of convergence of the p.p.e. on *regular splines*, when the Lévy density belongs to some *Lipschitz spaces* or *Besov spaces* of smooth functions, are derived. In Section 6, implementation of the method using discrete-time sampling of the process is briefly discussed. We finish with proofs of the main results.

## 2. A non-parametric estimation method

Consider a real Lévy process $X = \{X(t)\}_{t \geq 0}$ with Lévy density $p : \mathbb{R}_0 \to \mathbb{R}_+$, where $\mathbb{R}_0 \equiv \mathbb{R} \backslash \{0\}$. Then, $X$ is a càdlàg process with independent and stationary increments such that the characteristic function of its marginals is given by

$$(2.1) \quad \mathbb{E}\left[e^{iuX(t)}\right] = \exp\left\{t\left(iub - \frac{u^2\sigma^2}{2} + \int_{\mathbb{R}_0}\left\{e^{iux} - 1 - iux1_{[|x|\leq 1]}\right\}p(x)dx\right)\right\},$$

with $p : \mathbb{R}_0 \to \mathbb{R}_+$ such that

$$(2.2) \qquad\qquad\qquad \int_{\mathbb{R}_0}(1 \wedge x^2)p(x)dx < \infty.$$

Since $X$ is a càdlàg process, the set of its jump times

$$\left\{t > 0 : \Delta X(t) \equiv X(t) - X(t^-) \neq 0\right\}$$

is countable. Moreover, for Borel subsets $B$ of $[0, \infty) \times \mathbb{R}_0$,

$$(2.3) \qquad \mathcal{J}(B) \equiv \# \left\{ t > 0 : (t, X(t) - X(t^-)) \in B \right\},$$

is a well-defined random measure on $[0, \infty) \times \mathbb{R}_0$, where $\#$ denotes cardinality. The Lévy-Itô decomposition of the sample paths (see Theorem 19.2 of [27]) implies that $\mathcal{J}$ is a Poisson process on the Borel sets $\mathcal{B}([0, \infty) \times \mathbb{R}_0)$ with mean measure

$$(2.4) \qquad \mu(B) = \iint\limits_{B} p(x)\, dt\, dx.$$

Recall also that the stochastic integral of a deterministic function $f : \mathbb{R}_0 \to \mathbb{R}$ with respect to $\mathcal{J}$ is defined by

$$(2.5) \qquad I\,(f) \equiv \iint\limits_{[0,T] \times \mathbb{R}_0} f(x) \mathcal{J}(dt, dx) = \sum_{t \leq T} f(\Delta X(t)),$$

where this last expression is well defined if

$$\int_0^T \int_{\mathbb{R}_0} |f(x)| \mu(dt, dx) = T \int_{\mathbb{R}_0} |f(x)| p(x) dx < \infty;$$

see e.g. Chapter 10 in [19].

We consider the problem of estimating the Lévy density $p$ on a Borel set $D \in \mathcal{B}(\mathbb{R}_0)$ using a *projection estimation approach*. According to this paradigm, $p$ is estimated by estimating its best approximating function in a finite-dimensional linear space $\mathcal{S}$. The linear space $\mathcal{S}$ is taken so that it has good approximation properties for general classes of functions. Typical choices are piecewise polynomials or wavelets. Throughout, we make the following standing assumption.

**Assumption 1.** The Lévy measure $\nu(dx) \equiv p(x) dx$ is absolutely continuous with respect to a known measure $\eta$ on $\mathcal{B}(D)$ so that the Radon-Nikodym derivative

$$(2.6) \qquad \frac{d\nu}{d\eta}(x) = s(x), \;\; x \in D,$$

is positive, bounded, and satisfies

$$(2.7) \qquad \int_D s^2(x) \eta(dx) < \infty.$$

In that case, $s$ is called the Lévy density, on $D$, of the process with respect to the reference measure $\eta$.

**Remark 2.1.** Under the previous assumption, the measure $\mathcal{J}$ of (2.3), when restricted to $\mathcal{B}([0, \infty) \times D)$, is a Poisson process with mean measure

$$(2.8) \qquad \mu(B) = \iint\limits_{B} s(x)\, dt\, \eta(dx), \;\; B \in \mathcal{B}([0, \infty) \times D).$$

Our goal will be to estimate the Lévy density $s$, which itself could in turn be used to retrieve $p$ on $D$ via (2.6). To illustrate this strategy consider a continuous Lévy density $p$ such that

$$p(x) = O\left(x^{-1}\right), \;\; \text{as } x \to 0.$$

This type of densities satisfies the above assumption with respect to the measure $\eta(dx) = x^{-2}dx$ on domains of the form $D = \{x : 0 < |x| < b\}$. Clearly, an estimator $\hat{p}$ for the Lévy density $p$ can be generated from an estimator $\hat{s}$ for $s$ by fixing $\hat{p}(x) \equiv x^{-2}\hat{s}(x)$.

Let us now describe the main ingredients of our approach. Let $\mathcal{S}$ be a finite dimensional subspace of $L^2 \equiv L^2((D, \eta))$ equipped with the standard norm

$$\|f\|_\eta^2 \equiv \int_D f^2(x)\, \eta(dx).$$

The space $\mathcal{S}$ plays the role of an approximating linear model for the Lévy density $s$. Of course, under the $L^2$ norm, the best approximation of $s$ on $\mathcal{S}$ is the orthogonal projection defined by

$$(2.9) \qquad s^\perp(x) \equiv \sum_{i=1}^d \left( \int_D \varphi_i(y)s(y)\eta(dy) \right) \varphi_i(x),$$

where $\{\varphi_1, \ldots, \varphi_d\}$ is an arbitrary orthonormal basis of $\mathcal{S}$. The *projection estimator* of $s$ on $S$ is defined by

$$(2.10) \qquad \hat{s}(x) \equiv \sum_{i=1}^d \hat{\beta}_i \varphi_i(x),$$

where we fix

$$(2.11) \qquad \hat{\beta}_i \equiv \frac{1}{T} \iint\limits_{[0,T] \times D} \varphi_i(x)\mathcal{J}(dt, dx).$$

This is the most natural unbiased estimator for the *orthogonal projection $s^\perp$*. Notice also that $\hat{s}$ is independent of the specific orthonormal basis of $\mathcal{S}$. Indeed, the projection estimator is the unique solution to the minimization problem

$$\min_{f \in \mathcal{S}} \gamma_D(f),$$

where $\gamma_D : L^2((D, \eta)) \to \mathbb{R}$ is given by

$$(2.12) \qquad \gamma_D(f) \equiv -\frac{2}{T} \iint\limits_{[0,T] \times D} f(x)\, \mathcal{J}(dt, dx) + \int_D f^2(x)\, \eta(dx).$$

In the literature on model selection (see e.g. [7] and [25]), $\gamma_D$ is the so-called *contrast function*. The previous characterization also provides a mechanism to numerically evaluate $\hat{s}$ when an orthonormal basis of $\mathcal{S}$ is not explicitly available.

The following proposition provides both the first-order and the second-order properties of $\hat{s}$. These follow directly from the well-known formulas for the mean and variance of Poisson integrals (see e.g. [19] Chapter 10).

**Proposition 2.2.** *Under Assumption 1, $\hat{s}$ is an unbiased estimator for $s^\perp$ and its "mean-square error", defined by*

$$\chi^2 \equiv \|\hat{s} - s^\perp\|_\eta^2 = \int \left( \hat{s}(x) - s^\perp(x) \right)^2 \eta(dx),$$

*is such that*

$$\mathbb{E}\left[\chi^2\right] = \frac{1}{T}\sum_{i=1}^{d}\int_{D}\varphi_i^2(x)s(x)\,\eta(dx).$$
(2.13)

*The* risk *of $\hat{s}$ admits the decomposition*

$$\mathbb{E}\left[\|s - \hat{s}\|_\eta^2\right] = \|s - s^\perp\|_\eta^2 + \mathbb{E}\left[\chi^2\right].$$
(2.14)

The first term in (2.14), the *bias term*, accounts for the distance of the unknown function $s$ to the model $\mathcal{S}$, while the second term, the *variance term*, measures the error of estimation within the linear model $\mathcal{S}$. Notice that (2.13) is finite because $s$ is assumed bounded on $D$ and thus,

$$\mathbb{E}\left[\chi^2\right] \leq \|s\|_\infty \frac{d}{T}.$$
(2.15)

## 3. Model selection via penalized projection estimator

A crucial issue in the above approach is the selection of the approximating linear model $\mathcal{S}$. In principle, a "nice" density $s$ can be approximated closely by general linear models such as splines or wavelet. However, a more robust model $\mathcal{S}'$ containing $\mathcal{S}$ will result in a better approximation of $s$, but with a larger variance. This raises the natural problem of selecting one model, out of a collection of linear models $\{\mathcal{S}_m, m \in \mathcal{M}\}$, that approximately realizes the best trade-off between the risk of estimation within the model and the distance of the unknown Lévy density to the approximating model.

Let $\hat{s}_m$ and $s_m^\perp$ be respectively the projection estimator and the orthogonal projection of $s$ on $\mathcal{S}_m$. The following equation, readily derived from (2.14), gives insight on a sensible solution to the model selection problem:

$$\mathbb{E}\left[\|s - \hat{s}_m\|_\eta^2\right] = \|s\|_\eta^2 + \mathbb{E}\left[-\|\hat{s}_m\|_\eta^2 + \text{pen}(m)\right].$$
(3.1)

Here, $\text{pen}(m)$ is defined in terms of an orthonormal basis $\{\varphi_{1,m}, \ldots, \varphi_{d_m,m}\}$ of $\mathcal{S}_m$ by the equation:

$$\text{pen}(m) = \frac{2}{T^2}\iint\limits_{[0,T]\times D}\left(\sum_{i=1}^{d_m}\varphi_{i,m}^2(x)\right)\mathcal{J}(dt,dx).$$
(3.2)

Equation (3.1) shows that the risk of $\hat{s}_m$ moves "parallel" to the expectation of the *observable statistics* $-\|\hat{s}_m\|_\eta^2 + \text{pen}(m)$. This fact justifies to choose the model that minimizes such statistics. We will see later that other choices for $\text{pen}(\cdot)$ also give good results. Therefore, given a penalization function $\text{pen} : \mathcal{M} \to [0, \infty)$, we consider estimators of the form

$$\tilde{s} \equiv \hat{s}_{\hat{m}},$$
(3.3)

where $\hat{s}_m$ is the projection estimator on $\mathcal{S}_m$ and

$$\hat{m} \equiv \text{argmin}_{m\in\mathcal{M}}\left\{-\|\hat{s}_m\|_\eta^2 + \text{pen}(m)\right\}.$$

An estimator $\tilde{s}$ as in (3.3) is called a **penalized projection estimator** (p.p.e.) on the collection of linear models $\{\mathcal{S}_m, m \in \mathcal{M}\}$.

Methods of estimation based on the minimization of penalty functions have a long history in the literature on regression and density estimation (for instance, [1, 22], and [28]). The general idea is to choose, among a given collection of parametric models, the model that minimizes a loss function plus a penalty term that controls the variance term, which will forcefully increase as the approximating linear models become more detailed. Such penalized estimation was promoted for nonparametric density estimation in [8], and in the context of non-homogeneous Poisson processes in [25].

## 4. Risk bound and oracle inequalities

The penalization idea of the previous section provides a sensible criterion to select an estimator $\tilde{s} \equiv \hat{s}_{\hat{m}}$ out of the projection estimators $\{\hat{s}_m : m \in \mathcal{M}\}$ induced by a given collection of approximating linear models $\{\mathcal{S}_m, m \in \mathcal{M}\}$. Ideally, one wishes to choose that projection estimator $\hat{s}_{m^*}$ that minimizes the risk; namely, such that

$$(4.1) \qquad \mathbb{E}\left[\|s - \hat{s}_{m^*}\|_\eta^2\right] \leq \mathbb{E}\left[\|s - \hat{s}_m\|_\eta^2\right], \quad \text{for all } m \in \mathcal{M}.$$

Of course, to pick the best $\hat{s}_m$ is not feasible since $s$ is not available to actually compute and compare the risks. But, how bad would the risk of $\tilde{s}$ be compared to the best possible risk that can be achieved by projection estimators? One can aspire to achieve the smallest possible risk up to a constant. In other words, it is desirable that our estimator $\tilde{s}$ comply with an inequality of the form

$$(4.2) \qquad \mathbb{E}\left[\|s - \tilde{s}\|_\eta^2\right] \leq C \inf_{m \in \mathcal{M}} \mathbb{E}\left[\|s - \hat{s}_m\|_\eta^2\right],$$

for a constant $C$ independent of the linear models. The model $\mathcal{S}_{m^*}$ that achieves the minimal risk (using projection estimation) is the *oracle model* and inequalities of the type (4.2) are called *oracle inequalities*. Approximate oracle inequalities were proved in [25] for the intensity function of a nonhomogeneous Poisson process. In this section we show that for certain penalization functions, the resulting penalized projection estimator $\tilde{s}$ defined by (3.3) satisfies the inequality

$$(4.3) \qquad \mathbb{E}\left[\|s - \tilde{s}\|_\eta^2\right] \leq C \inf_{m \in \mathcal{M}} \mathbb{E}\left[\|s - \hat{s}_m\|_\eta^2\right] + \frac{C'}{T},$$

for some "model free" constants $C, C'$ (remember that the time period of observations is $[0, T]$). The main tool in obtaining oracle inequalities is an upper bound for the risk of the penalized projection estimator $\tilde{s}$. The proof of (4.3) follows essentially from the arguments in [25]; however, to overcome the possible lack of finiteness on the reference measure $\eta$ (see Assumption 1), which is required in [25], and to avoid superfluous rough upper bounds, the dimension of the linear model is explicitly included in the penalization and the arguments are refined.

Let us introduce some notation. Below, $d_m$ denotes the dimension of the linear model $\mathcal{S}_m$, and $\{\varphi_{1,m}, \ldots, \varphi_{d_m,m}\}$ is an arbitrary orthonormal basis of $\mathcal{S}_m$. Define

$$(4.4) \qquad D_m = \sup\left\{\|f\|_\infty^2 : f \in S_m, \|f\|_\eta^2 = 1\right\},$$

which is assumed to be finite and can be proved to be equal to $\|\sum_{i=1}^{d_m} \varphi_{i,m}^2\|_\infty$.

We make the following regularity condition, introduced in [25], that essentially controls the complexity of the linear models. This assumption is satisfied by splines and trigonometric polynomials, but not by wavelet bases.

**Assumption 2.** There exist constants $\Gamma > 0$ and $R \geq 0$ such that for every positive integer $n$,

$$\# \{ m \in \mathcal{M} : d_m = n \} \leq \Gamma n^R.$$

We now present our main result.

**Theorem 4.1.** *Let* $\{ \mathcal{S}_m, m \in \mathcal{M} \}$ *be a family of finite dimensional linear subspaces of* $L^2((D, \eta))$ *satisfying Assumption 2 and such that* $D_m < \infty$. *Let* $\mathcal{M}_T \equiv \{ m \in \mathcal{M} : D_m \leq T \}$. *If* $\hat{s}_m$ *and* $s_m^\perp$ *are respectively the projection estimator and the orthogonal projection of the Lévy density* $s$ *on* $\mathcal{S}_m$ *then, the penalized projection estimator* $\tilde{s}_T$ *on* $\{ \mathcal{S}_m \}_{m \in \mathcal{M}_T}$ *defined by* (3.3) *is such that*

$$(4.5) \qquad \mathbb{E}\left[ \| s - \tilde{s}_T \|_\eta^2 \right] \leq C \inf_{m \in \mathcal{M}_T} \left\{ \| s - s_m^\perp \|_\eta^2 + \mathbb{E}\left[ \mathrm{pen}(m) \right] \right\} + \frac{C'}{T},$$

*whenever* pen $: \mathcal{M} \to [0, \infty)$ *takes either one of the following forms for some fixed (but arbitrary) constants* $c > 1$, $c' > 0$, *and* $c'' > 0$:
**(a)** $\mathrm{pen}(m) \geq c \frac{D_m \mathcal{N}}{T^2} + c' \frac{d_m}{T}$, *where* $\mathcal{N} \equiv \mathcal{J}([0, T] \times D)$ *is the number of jumps prior to* $T$ *with sizes in* $D$ *and where it is assumed that* $\rho \equiv \int_D s(x) \eta(dx) < \infty$;
**(b)** $\mathrm{pen}(m) \geq c \frac{\hat{V}_m}{T}$, *where* $\hat{V}_m$ *is defined by*

$$(4.6) \qquad \hat{V}_m \equiv \frac{1}{T} \iint\limits_{[0,T] \times D} \left( \sum_{i=1}^{d_m} \varphi_{i,m}^2(x) \right) \mathcal{J}(dt, dx),$$

*and where it is assumed that* $\beta \equiv \inf_{m \in \mathcal{M}} \frac{\mathbb{E}[\hat{V}_m]}{D_m} > 0$ *and that* $\phi \equiv \inf_{m \in \mathcal{M}} \frac{D_m}{d_m} > 0$;
**(c)** $\mathrm{pen}(m) \geq c \frac{\hat{V}_m}{T} + c' \frac{D_m}{T} + c'' \frac{d_m}{T}$.
In (4.5), *the constant* $C$ *depends only on* $c$, $c'$ *and* $c''$, *while* $C'$ *varies with* $c$, $c'$, $c''$, $\Gamma$, $R$, $\| s \|_\eta$, $\| s \|_\infty$, $\rho$, $\beta$, *and* $\phi$.

**Remark 4.2.** It can be shown that if $c \geq 2$, then for arbitrary $\varepsilon > 0$, there is a constant $C'(\varepsilon)$ (increasing as $\varepsilon \downarrow 0$) such that

$$(4.7) \qquad \mathbb{E}\| s - \tilde{s} \|_\eta^2 \leq (1 + \varepsilon) \inf_{m \in \mathcal{M}} \left\{ \| s - s_m^\perp \|_\eta^2 + \mathbb{E}\left[ \mathrm{pen}(m) \right] \right\} + \frac{C'(\varepsilon)}{T}.$$

One important consequence of the risk bound (4.5) is the following oracle inequality:

**Corollary 4.3.** *In the setting of Theorem* 4.1(b), *if the penalty function is of the form* $\mathrm{pen}(m) \equiv c \frac{\hat{V}_m}{T}$, *for every* $m \in \mathcal{M}_T$, $\beta > 0$, *and* $\phi > 0$, *then*

$$(4.8) \qquad \mathbb{E}\left[ \| s - \tilde{s}_T \|_\eta^2 \right] \leq \widetilde{C} \inf_{m \in \mathcal{M}_T} \left\{ \mathbb{E}\left[ \| s - \hat{s}_m \|_\eta^2 \right] \right\} + \frac{\widetilde{C}'}{T},$$

*for a constant* $\widetilde{C}$ *depending only on* $c$, *and a constant* $\widetilde{C}'$ *depending on* $c$, $\Gamma$, $R$, $\| s \|_\eta$, $\| s \|_\infty$, $\beta$, *and* $\phi$.

## 5. Rate of convergence for smooth Lévy densities

We use the risk bound of the previous section to study the "long run" ($T \to \infty$) rate of convergence of penalized projection estimators based on regular piecewise polynomials, when the Lévy density is "smooth". More precisely, on a window of estimation $D \equiv [a,b] \subset \mathbb{R}_0$, the Lévy density of the process with respect to the Lebesgue measure $\eta(dx) \equiv dx$, denoted by $s$, is assumed to belong to the *Besov space* (also called *Lipschitz* space) $\mathcal{B}_\infty^\alpha (L^p([a,b]))$ for some $p \in [2, \infty]$ and $\alpha > 0$ (see for instance [15] and references therein for background on these spaces). Concretely, $\mathcal{B}_\infty^\alpha (L^p([a,b]))$ consists of those functions $f \in L^p([a,b], dx)$ if $0 < p < \infty$ (or $f$ continuous if $p = \infty$) such that

$$|f|_{\mathcal{B}_\infty^\alpha(L^p)} \equiv \sup_{\delta > 0} \frac{1}{\delta^\alpha} \sup_{0 < h \leq \delta} \|\Delta_h^r(f, \cdot)\|_{L^p([a,b],dx)} < \infty,$$

where $\Delta_h(f, x) \equiv f(x+h) - f(x)$ and $\Delta_h^r(f, x)$ is the $r^{th}$-order difference of $f$ defined by

$$\Delta_h^r(f, x) \equiv \Delta_h(\Delta_h^{r-1}(f, \cdot), x),$$

for $x$ such that $x + rh \in D$ and $r \in \mathbb{N}$. The following spaces are closely related. For $k \in \mathbb{N}$ and $\beta \in (0, 1]$ such that $\alpha = k + \beta$, let $\text{Lip}(\alpha, L^p([a,b]))$ be the class of functions $f$ such that $f, \ldots, f^{(k-1)}$ are absolutely continuous on $[a,b]$ with $f^{(k)} \in L^p((a,b))$ satisfying

$$\|\Delta_h(f^{(k)}, \cdot)\|_{L^p([a,b],dx)} \leq M h^\beta,$$

for some $M < \infty$. It is know that if $\alpha > 0$ is not an integer and $1 \leq p \leq \infty$, then $f \in \text{Lip}(\alpha, L^p([a,b]))$ if and only if $f$ is a.e. equal to a function in $\mathcal{B}_\infty^\alpha (L^p([a,b]))$. In general, $\text{Lip}(\alpha, L^p([a,b])) \subset \mathcal{B}_\infty^\alpha (L^p([a,b]))$, for any $0 < p \leq \infty$ and $\alpha > 0$ (see e.g. [15]).

An important reason for the choice of the Besov class of smooth functions is the availability of estimates for the error of approximation by splines, trigonometric polynomials, and wavelets (see e.g. [15] and [4]). In particular, if $\mathcal{S}_m^k$ denotes the space of piecewise polynomials of degree at most $k$, based on the regular partition of $[a,b]$ with $m$ intervals ($m \geq 1$), and $s \in \mathcal{B}_\infty^\alpha (L^p([a,b]))$ with $k > \alpha - 1$, then there exists a constant $C(s)$ such that

$$(5.1) \qquad \mathrm{d}_p\left(s, \mathcal{S}_m^k\right) \leq C(s) m^{-\alpha},$$

where $\mathrm{d}_p$ is the distance induced by the $L^p$-norm on $([a,b], dx)$ (see [15]). The following gives the rate of convergence of the p.p.e. on regular splines.

**Corollary 5.1.** *With the notation of Theorem 4.1, taking $D = [a,b]$ and $\eta(dx) = dx$ , let $\tilde{s}_T$ be the penalized projection estimator on $\{\mathcal{S}_m^k\}_{m \in \mathcal{M}_T}$ with penalization*

$$\text{pen}(m) \equiv c \frac{\hat{V}_m}{T} + c' \frac{D_m}{T} + c'' \frac{d_m}{T},$$

*for some fixed $c > 1$ and $c', c'' > 0$. Then, if the restriction to $D$ of the Lévy density $s$ belongs to $\mathcal{B}_\infty^\alpha (L^p([a,b]))$, with $2 \leq p \leq \infty$ and $0 < \alpha < k + 1$, then*

$$\limsup_{T \to \infty} T^{2\alpha/(2\alpha+1)} \mathbb{E}\left[\|s - \tilde{s}_T\|_\eta^2\right] < \infty.$$

*Moreover, for any $R > 0$ and $L > 0$,*

$$(5.2) \qquad \limsup_{T \to \infty} T^{2\alpha/(2\alpha+1)} \sup_{s \in \Theta(R,L)} \mathbb{E}\left[\|s - \tilde{s}_T\|_\eta^2\right] < \infty,$$

where $\Theta(R, L)$ consists of all the Lévy densities $f$ such that $\|f\|_{L^\infty([a,b],dx)} < R$, and such that the restriction of $f$ to $[a,b]$ is a member of $\mathcal{B}^\alpha_\infty (L^p([a,b]))$ with $|f|_{\mathcal{B}^\alpha_\infty(L^p)} < L$.

The previous result implies that the p.p.e. on regular splines has a rate of convergence of order $T^{-2\alpha/(2\alpha+1)}$ for the class of Besov Lévy densities $\Theta(R, L)$.

## 6. Estimation based on discrete time data

Let us finish with some remarks on how to approximate the continuous-time statistics of our methods using only discrete-time observations. In practice, we can aspire to sample the process $X(t)$ at discrete times, but we are neither able to measure the size of the jumps $\Delta X(t) \equiv X(t) - X(t^-)$ nor the times of the jumps $\{t : \Delta X(t) > 0\}$. In general, Poisson integrals of the type

$$(6.1) \qquad I(f) \equiv \iint_{[0,T] \times \mathbb{R}_0} f(x) \mathcal{J}(dt, dx) = \sum_{t \le T} f(\Delta X(t)),$$

are not accessible. Intuitively, the following statistic is the most natural approximation to (6.1):

$$(6.2) \qquad I_n(f) \equiv \sum_{k=1}^n f(\Delta_k X),$$

where $\Delta_k X$ is the $k^{th}$ increment of the process with time span $h_n \equiv T/n$; that is,

$$\Delta_k X \equiv X(kh_n) - X((k-1)h_n), \quad k = 1, \ldots, n.$$

How good is this approximation and in what sense? Under some conditions on $f$, we can readily prove the weak convergence of (6.2) to (6.1) using properties of the transition distributions of $X$ in small time (see [5], Corollary 8.9 of [27], and [26]). The following theorem summarizes some known results on the small-time transition distribution.

**Theorem 6.1.** *Let $X = \{X(t)\}_{t \ge 0}$ be a Lévy process with Lévy measure $\nu$. The following statements hold true.*

**(1)** *For each $a > 0$,*

$$(6.3) \qquad \lim_{t \to 0} \frac{1}{t} \mathbb{P}(X(t) > a) = \nu([a, \infty)),$$

$$(6.4) \qquad \lim_{t \to 0} \frac{1}{t} \mathbb{P}(X(t) \le -a) = \nu((-\infty, -a]).$$

**(2)** *For any continuous bounded function $h$ vanishing in a neighborhood of the origin,*

$$(6.5) \qquad \lim_{t \to 0} \frac{1}{t} \mathbb{E}[h(X(t))] = \int_{\mathbb{R}_0} h(x) \nu(dx).$$

**(3)** *If $h$ is continuous and bounded and if $\lim_{|x| \to 0} h(x)|x|^{-2} = 0$, then*

$$\lim_{t \to 0} \frac{1}{t} \mathbb{E}[h(X(t))] = \int_{\mathbb{R}_0} h(x) \nu(dx).$$

*Moreover, if $\int_{\mathbb{R}_0} (|x| \wedge 1) \nu(dx) < \infty$, it suffices to have $h(x)(|x| \wedge 1)^{-1}$ continuous and bounded.*

Convergence results like (6.5) are useful to establish the convergence in distribution of $I_n(f)$ since

$$\mathbb{E}\left[e^{\mathrm{i}uI_n(f)}\right] = \left(\mathbb{E}\left[e^{\mathrm{i}uf\left(X\left(\frac{T}{n}\right)\right)}\right]\right)^n = \left(1 + \frac{a_n}{n}\right)^n,$$

where $a_n = n\mathbb{E}\left[h\left(X\left(\frac{T}{n}\right)\right)\right]$ with $h(x) = e^{\mathrm{i}uf(x)} - 1$. So, if $f$ is such that

$$(6.6) \qquad \lim_{t\to 0}\frac{1}{t}\mathbb{E}\left[e^{\mathrm{i}uf(X(t))} - 1\right] = \int_{\mathbb{R}_0}\left(e^{\mathrm{i}uf(x)} - 1\right)\nu(dx),$$

then $a_n$ converges to $a \equiv T\int_{\mathbb{R}_0}h(x)\nu(dx)$, and thus

$$\lim_{n\to\infty}\left(1 + \frac{a_n}{n}\right)^n = \lim_{n\to\infty}e^{n\log\left(1+\frac{a_n}{n}\right)} = e^a.$$

We thus have the following result.

**Proposition 6.2.** *Let* $X = \{X(t)\}_{t\geq 0}$ *be a Lévy process with Lévy measure* $\nu$. *Then,*

$$\lim_{n\to\infty}\mathbb{E}\left[e^{\mathrm{i}uI_n(f)}\right] = \exp\left\{T\int_{\mathbb{R}_0}\left(e^{\mathrm{i}uf(x)} - 1\right)\nu(dx)\right\},$$

*if* $f$ *satisfies either one of the following conditions:*

(1) $f(x) = \mathbf{1}_{(a,b]}(x)h(x)$ *for an interval* $[a,b] \subset \mathbb{R}_0$ *and a continuous function* $h$;
(2) $f$ *is continuous on* $\mathbb{R}_0$ *and* $\lim_{|x|\to 0}f(x)|x|^{-2} = 0$.

*In particular,* $I_n(f)$ *converges in distribution to* $I(f)$ *under any one of the previous two conditions.*

**Remark 6.3.** Notice that if (6.5) holds true when replacing $h$ by $f$ and $f^2$, then the mean and variance of $I_n(f)$ obey the asymptotics:

$$\lim_{n\to\infty}\mathbb{E}\left[I_n(f)\right] = T\int_{\mathbb{R}_0}f(x)\nu(dx);$$

$$\lim_{n\to\infty}\mathrm{Var}\left[I_n(f)\right] = T\int_{\mathbb{R}_0}f^2(x)\nu(dx).$$

**Remark 6.4.** Very recently, [23] proposed a procedure to disentangle the jumps from the diffusion part in the case of jump-diffusion models driven by finite-jump activity Lévy processes. It is proved there that for certain functions $r : \mathbb{R}_+ \to \mathbb{R}_+$, there exists $N(\omega)$ such that for $n \geq N(\omega)$, a jump occurs in the interval $((k-1)h_n, kh_n]$ if and only if $(\Delta_k X)^2 > r(h_n)$. Here, $h_n = T/n$ and $\Delta_k X$ is the $k^{th}$ increment of the process. These results suggest to use statistics of the form

$$\sum_{k=1}^n f\left(\Delta_k X\right)\mathbf{1}\left[\left(\Delta_k X\right)^2 > r(h_n)\right]$$

instead of (6.2) to approximate the integral (6.1).

## 7. Proofs

### 7.1. Proof of the risk bound

We break the proof of Theorem 4.1 into several preliminary results.

**Lemma 7.1.** *For any penalty function* pen : $\mathcal{M} \to [0, \infty)$ *and any* $m \in \mathcal{M}$, *the penalized projection estimator* $\tilde{s}$ *satisfies*

$$(7.1) \qquad \|s - \tilde{s}\|_\eta^2 \le \|s - s_m^\perp\|_\eta^2 + 2\chi_{\hat{m}}^2 + 2\nu_D \left(s_{\hat{m}}^\perp - s_m^\perp\right) + \text{pen}(m) - \text{pen}(\hat{m}),$$

*where* $\chi_m^2 \equiv \|s_m^\perp - \hat{s}_m\|_\eta^2$ *and where the functional* $\nu_D : L^2\left((D, \eta)\right) \to \mathbb{R}$ *is defined by*

$$(7.2) \qquad \nu_D(f) \equiv \iint\limits_{[0,T] \times D} f(x) \; \frac{\mathcal{J}(dt, dx) - s(x)\, dt\, \eta(dx)}{T}.$$

The general idea in obtaining (4.5) is to bound the "inaccessible" terms on the right hand side of (7.1) (namely $\chi_{\hat{m}}^2$ and $\nu_D \left(s_{\hat{m}}^\perp - s_m^\perp\right)$) by observable statistics. In fact, the penalizations pen($\cdot$) given in Theorem 4.1 are chosen so that the right hand side in (7.1) does not involve $\hat{m}$. To carry out this plan, we use concentration inequalities for $\chi_{\hat{m}}^2$ and for the compensated Poisson integrals $\nu_D(f)$. The following result gives a concentration inequality for general compensated Poisson integrals.

**Proposition 7.2.** *Let $N$ be a Poisson process on a measurable space* (V, $\mathcal{V}$) *with mean measure $\mu$ and let $f :$ V $\to \mathbb{R}$ be an essentially bounded measurable function satisfying* $0 < \|f\|_\mu^2 \equiv \int_V f^2(v)\mu(dv)$ *and* $\int_V |f(v)|\mu(dv) < \infty$. *Then, for any* $u > 0$,

$$(7.3) \qquad \mathbb{P}\left[\int_V f(v)(N(dv) - \mu(dv)) \ge \|f\|_\mu \sqrt{2u} + \frac{1}{3}\|f\|_\infty u\right] \le e^{-u}.$$

*In particular, if $f :$ V $\to [0, \infty)$ then, for any $\epsilon > 0$ and $u > 0$,*
$$(7.4)$$
$$\mathbb{P}\left[(1 + \varepsilon)\left(\int_V f(v)N(dv) + \left(\frac{1}{2\varepsilon} + \frac{5}{6}\right)\|f\|_\infty u\right) \ge \int_V f(v)\mu(dv)\right] \ge 1 - e^{-u}.$$

For a proof of the inequality (7.3), see [25] (Proposition 7) or [18] (Corollary 5.1). Inequality (7.4) is a direct consequence of (7.3) (see Section 7.2 for a proof).

The next result allows us to bound the Poisson functional $\chi_m^2$. This result is essentially Proposition 9 of [25].

**Lemma 7.3.** *Let $N$ be a Poisson process on a measurable space* (V, $\mathcal{V}$) *with mean measure $\mu(dv) = p(v)\zeta(dv)$ and intensity function $p \in L^2(V, \mathcal{V}, \zeta)$. Let $\mathcal{S}$ be a finite dimensional subspace of $L^2(V, \mathcal{V}, \zeta)$ with orthonormal basis $\{\tilde{\varphi}_1, \ldots, \tilde{\varphi}_d\}$, and let*

$$(7.5) \qquad \hat{p}(v) \equiv \sum_{i=1}^d \left(\int_V \tilde{\varphi}_i(w)N(dw)\right) \tilde{\varphi}_i(v)$$

$$(7.6) \qquad p^\perp(v) \equiv \sum_{i=1}^d \left(\int_V p(w)\tilde{\varphi}_i(w)\eta(dw)\right) \tilde{\varphi}_i(v).$$

*Then, $\chi^2(\mathcal{S}) \equiv \|\hat{p} - p^\perp\|_\zeta^2$ is such that for any $u > 0$ and $\varepsilon > 0$*

$$(7.7) \qquad \mathbb{P}\left[\chi(\mathcal{S}) \ge (1 + \varepsilon)\sqrt{\mathbb{E}\left[\chi^2(\mathcal{S})\right]} + \sqrt{2kM_\mathcal{S}u} + k(\varepsilon)B_\mathcal{S}u\right] \le e^{-u},$$

*where we can take $k = 6$, $k(\varepsilon) = 1.25 + 32/\varepsilon$, and where*

$$(7.8) \qquad M_{\mathcal{S}} \equiv \sup \left\{ \int_V f^2(v) p(v) \zeta(dv) : f \in \mathcal{S}, \|f\|_\varsigma = 1 \right\},$$

$$(7.9) \qquad B_{\mathcal{S}} \equiv \sup \left\{ \|f\|_\infty : f \in \mathcal{S}, \|f\|_\varsigma = 1 \right\}.$$

Following the same strategy as in [25], the idea is to obtain from the previous lemmas a concentration inequality of the form

$$\mathbb{P}\left[ \|s - \tilde{s}\|_\eta^2 \le C \left( \|s - s_m^\perp\|_\eta^2 + \text{pen}(m) \right) + h(\xi) \right] \ge 1 - C' e^{-\xi},$$

for constants $C$ and $C'$, and a function $h(\xi)$ (all independent of $m$). This will prove to be enough in view of the following elementary result (see Section 7.2 for a proof).

**Lemma 7.4.** *Let $h : [0, \infty) \to \mathbb{R}_+$ be a strictly increasing function with continuous derivative such that $h(0) = 0$ and $\lim_{\xi \to \infty} e^{-\xi} h(\xi) = 0$. If $Z$ is random variable satisfying*

$$\mathbb{P}\left[ Z \ge h(\xi) \right] \le K e^{-\xi},$$

*for every $\xi > 0$, then*

$$\mathbb{E} Z \le K \int_0^\infty e^{-u} h(u) du.$$

We are now in a position to prove Theorem 4.1. Throughout the proof, we will have to introduce various constants and inequalities that will hold with high probability. In order to clarify the role that the constants play in these inequalities, we shall make some convention and give to the letters $x$, $y$, $f$, $a$, $b$, $\xi$, $\mathcal{K}$, c, and $C$, with various sub- or superscripts, special meaning. The letters with $x$ are reserved to denote positive constants that can be chosen arbitrarily. The letters with $y$ denote arbitrary constants greater than 1. $f, f_1, f_2, \ldots$ denote quadratic polynomials of the variable $\xi$ whose coefficients (denoted by $a's$ and $b's$) are determined by the values of the $x's$ and $y's$. The inequalities will be true with probabilities greater that $1 - \mathcal{K} e^{-\xi}$, where $\mathcal{K}$ is determined by the values of the $x's$ and the $y's$. Finally, $c's$ and $C's$ are used to denote constants constrained by the $x's$ and $y's$. *It is important to remember that the constants in a given inequality are meant only for that inequality.* The pair of equivalent inequalities below will be repeatedly invoked throughout the proof:

$$(7.10) \qquad \begin{array}{ll} \text{(i)} & 2ab \le xa^2 + \frac{1}{x}b^2, \quad \text{and} \\ \text{(ii)} & (a+b)^2 \le (1+x) a^2 + \left(1 + \frac{1}{x}\right) b^2, \quad \text{(for } x > 0). \end{array}$$

*Also, for simplicity, we write below $\|\cdot\|$ to denote the $L^2-$norm with respect to the reference measure $\eta$.*

**Proof of Theorem 4.1**. We consider successive improvements of the inequality (7.1):

**Inequality 1.** *For any positive constants $x_1$, $x_2$, $x_3$, and $x_4$, there exist a positive number $\mathcal{K}$ and an increasing quadratic function $f$ (both independent of the family of linear models and of $T$) such that, with probability larger than $1 - \mathcal{K} e^{-\xi}$,*

$$(7.11) \qquad \begin{aligned} \|s - \tilde{s}\|^2 &\le \|s - s_m^\perp\|^2 + 2\chi_{\hat{m}}^2 + 2x_1 \|s_{\hat{m}}^\perp - s_m^\perp\|^2 \\ &\quad + x_2 \frac{D_{\hat{m}}}{T} + x_3 \frac{D_m}{T} + x_4 \frac{d_{\hat{m}}}{T} \\ &\quad + \text{pen}(m) - \text{pen}(\hat{m}) + \frac{f(\xi)}{T}. \end{aligned}$$

*Proof.* Let us find an upper bound for $\nu_D(s_{m'}^\perp - s_m^\perp)$, $m', m \in \mathcal{M}$. Since the operator $\nu_D$ defined by (7.2) is just a compensated integral with respect to a Poisson process with mean measure $\mu(dtdx) = dt\eta(dx)$, we can apply Proposition 7.2 to obtain that, for any $x_{m'}' > 0$, and with probability larger than $1 - e^{-x_{m'}'}$

$$(7.12) \qquad \nu_D\left(s_{m'}^\perp - s_m^\perp\right) \leq \left\|\frac{s_{m'}^\perp - s_m^\perp}{T}\right\|_\mu \sqrt{2x_{m'}'} + \frac{\|s_{m'}^\perp - s_m^\perp\|_\infty x_{m'}'}{3T}.$$

In that case, the probability that (7.12) holds for every $m' \in \mathcal{M}$ is larger than $1 - \sum_{m' \in \mathcal{M}} e^{-x_{m'}}$ because $P(A \cap B) \geq 1 - a - b$, whenever $P(A) \geq 1 - a$ and $P(B) \geq 1 - b$. Clearly,

$$\left\|\frac{s_{m'}^\perp - s_m^\perp}{T}\right\|_\mu^2 = \iint_{[0,T]\times D} \left(\frac{s_{m'}^\perp(x) - s_m^\perp(x)}{T}\right)^2 s(x)dt\eta(dx)$$

$$\leq \|s\|_\infty \frac{\|s_{m'}^\perp - s_m^\perp\|^2}{T}.$$

Using (7.10)(i), the first term on the right hand side of (7.12) is then bounded as follows:

$$(7.13) \qquad \left\|\frac{s_{m'}^\perp - s_m^\perp}{T}\right\|_\mu \sqrt{2x_{m'}'} \leq x_1\|s_{m'}^\perp - s_m^\perp\|^2 + \frac{\|s\|_\infty x_{m'}'}{2Tx_1},$$

for any $x_1 > 0$. Using (4.4) and (7.10-i),

$$\|s_{m'}^\perp - s_m^\perp\|_\infty x_{m'}' \leq \left(\|s_{m'}^\perp\|_\infty + \|s_m^\perp\|_\infty\right) x_{m'}'$$

$$\leq \left(\sqrt{D_{m'}}\|s_{m'}^\perp\| + \sqrt{D_m}\|s_m^\perp\|\right) x_{m'}'$$

$$\leq \sqrt{D_{m'}}\|s\|x_{m'}' + \sqrt{D_m}\|s\|x_{m'}'$$

$$\leq 3x_2 D_{m'} + 3x_3 D_m + \frac{\|s\|^2 x_{m'}'^2}{12}\left(\frac{1}{x_2} + \frac{1}{x_3}\right),$$

for all $x_2 > 0$, $x_3 > 0$. It follows that, for any $x_1 > 0$, $x_2 > 0$, and $x_3 > 0$,

$$\nu_D\left(s_{m'}^\perp - s_m^\perp\right) \leq x_1\|s_{m'}^\perp - s_m^\perp\|^2 + x_2\frac{D_{m'}}{T} + x_3\frac{D_m}{T}$$

$$+ \frac{\|s\|_\infty x_{m'}'}{2Tx_1} + \frac{\|s\|^2 x_{m'}'^2}{36T\bar{x}},$$

where we set $\frac{1}{\bar{x}} = \frac{1}{x_2} + \frac{1}{x_3}$. Next, take

$$x_{m'}' \equiv x_4\sqrt{d_{m'}}\left(\frac{1}{\|s\|} \wedge \frac{1}{\|s\|_\infty}\right) + \xi.$$

Then, for any positive $x_1$, $x_2$, $x_3$, and $x_4$, there is a $\mathcal{K}$ and a function $f$ such that, with probability greater than $1 - \mathcal{K}e^{-\xi}$,

$$(7.14) \qquad \nu_D\left(s_{m'}^\perp - s_m^\perp\right) \leq x_1\|s_{m'}^\perp - s_m^\perp\|^2 + x_2\frac{D_{m'}}{T} + x_3\frac{D_m}{T}$$

$$+ \left(\frac{x_4^2}{18\bar{x}} + \frac{x_4}{2x_1}\right)\frac{d_{m'}}{T} + \frac{f(\xi)}{T}, \quad \forall m' \in \mathcal{M}.$$

Concretely,

$$(7.15)\qquad \begin{aligned} f(\xi) &= \frac{\|s\|}{18\bar{x}}\xi^2 + \frac{\|s\|_\infty}{2x_1}\xi, \\ \mathcal{K} &= \Gamma \sum_{n=1}^\infty n^R \exp\left(-\sqrt{n}x_4\left(\frac{1}{\|s\|} \wedge \frac{1}{\|s\|_\infty}\right)\right). \end{aligned}$$

Here, we used the assumption of polynomial models (Definition 2) to come up with the constant $\mathcal{K}$. Plugging (7.14) in (7.1), and renaming the coefficient of $d_{m'}/T$, we can corroborate inequality 1. $\qquad \square$

**Inequality 2.** *For any positive constants $y_1 > 1$, $x_2$, $x_3$, and $x_4$, there are positive constants $C_1 < 1$, $C_1' > 1$, and $\mathcal{K}$, and a strictly increasing quadratic polynomial $f$ (all independent of the class of linear models and of $T$) such that with probability larger than $1 - \mathcal{K}e^{-\xi}$,*

$$(7.16)\qquad \begin{aligned} C_1\|s-\tilde{s}\|^2 \leq\ &C_1'\|s-s_m^\perp\|^2 + y_1\chi_{\hat{m}}^2 \\ &+ x_2\frac{D_{\hat{m}}}{T} + x_3\frac{D_m}{T} + x_4\frac{d_{\hat{m}}}{T} \\ &+ \operatorname{pen}(m) - \operatorname{pen}(\hat{m}) + \frac{f(\xi)}{T}. \end{aligned}$$

*Moreover, if $1 < y_1 < 2$, then $C_1' = 3 - y_1$ and $C_1 = y_1 - 1$. If $y_1 \geq 2$, then $C_1' = 1 + 4x_1$ and $C_1 = 1 - 4x_1$, where $x_1$ is any positive constant related to $f$ via to the equation (7.15).*

*Proof.* Let us combine the term on the left hand side of (7.11) with the first three terms on the right hand side. Using the triangle inequality followed by (7.10-ii),

$$\|s_{\hat{m}}^\perp - s_m^\perp\|^2 \leq 2\|s-s_m^\perp\|^2 + 2\|s_{\hat{m}}^\perp - s\|^2.$$

Then, since $\chi_{\hat{m}}^2 = \|s_{\hat{m}}^\perp - \hat{s}_{\hat{m}}\|^2$, and $\|s_{\hat{m}}^\perp - s\|^2 = \|s - \hat{s}_{\hat{m}}\|^2 - \|s_{\hat{m}}^\perp - \hat{s}_{\hat{m}}\|^2$, it follows that

$$\begin{aligned} \|s - s_m^\perp\|^2 &+ 2\chi_{\hat{m}}^2 + 2x_1\|s_{\hat{m}}^\perp - s_m^\perp\|^2 - \|s-\tilde{s}\|^2 \\ &\leq (1+4x_1)\|s-s_m^\perp\|^2 + (2-4x_1)\|s_{\hat{m}}^\perp - \hat{s}_{\hat{m}}\|^2 \\ &\quad + (4x_1 - 1)\|s-\tilde{s}\|^2, \end{aligned}$$

for every $x_1 > 0$. Then, for any $y_1 > 1$, there are positive constants $C > 0$, $C_1' > 1$, and $C_1 < 1$ such that

$$(7.17)\qquad \begin{aligned} \|s - s_m^\perp\|^2 &+ 2\chi_{\hat{m}}^2 + 2C\|s_{\hat{m}}^\perp - s_m^\perp\|^2 - \|s-\tilde{s}\|^2 \\ &\leq C_1'\|s-s_m^\perp\|^2 + y_1\chi_{\hat{m}}^2 - C_1\|s-\tilde{s}\|^2. \end{aligned}$$

Combining (7.11) and (7.17), we obtain (7.16). $\qquad \square$

**Inequality 3.** *For any $y_2 > 1$ and positive constants $x_i$, $i = 2, 3, 4$, there exist positive reals $C_1 < 1$, $C_1' > 1$, an increasing quadratic polynomial of the form $f_2(\xi) = a\xi^2 + b\xi$, and a constant $\mathcal{K}_2 > 0$ (all independent of the family of linear models and of $T$) so that, with probability greater than $1 - \mathcal{K}_2e^{-\xi}$,*

$$(7.18)\qquad \begin{aligned} C_1\|s-\tilde{s}\|^2 \leq\ &C_1'\|s-s_m^\perp\|^2 \\ &+ y_2\frac{V_{\hat{m}}}{T} + x_2\frac{D_{\hat{m}}}{T} + x_3\frac{d_{\hat{m}}}{T} - \operatorname{pen}(\hat{m}) \\ &+ x_4\frac{D_m}{T} + \operatorname{pen}(m) + \frac{f(\xi)}{T}. \end{aligned}$$

*Proof.* We bound $\chi^2_{m'}$ using Lemma 7.3 with $V = \mathbb{R}_+ \times D$ and $\mu(d\mathbf{x}) = s(x)dt\eta(dx)$. We regard the linear model $\mathcal{S}_m$ as a subspace of $L^2(\mathbb{R}_+ \times D, dt\eta(dx))$ with orthonormal basis $\{\frac{\varphi_{1,m}}{\sqrt{T}}, \ldots, \frac{\varphi_{d_m,m}}{\sqrt{T}}\}$. Recall that

$$\chi^2_m = \|s^\perp_m - \hat{s}_m\|^2 = \sum_{i=1}^d \left[ \iint_{[0,T] \times D} \varphi_{i,m}(x) \frac{\mathcal{J}(dt, dx) - s(x)dt\eta(dx)}{T} \right]^2.$$

Then, with probability larger than $1 - \sum_{m' \in \mathcal{M}} e^{-x'_{m'}}$,

(7.19) $$\sqrt{T}\chi_{m'} \leq (1 + x_1)\sqrt{V_{m'}} + \sqrt{2kM_{m'}x'_{m'}} + k(x_1)B_{m'}x'_{m'},$$

for every $m' \in \mathcal{M}$, where $B_{m'} = \sqrt{D_{m'}/T}$,

(7.20) $$V_{m'} \equiv \int_D \left( \sum_{i=1}^{d_m} \varphi^2_{i,m}(x) \right) s(x)\eta(dx),$$

$$M_{m'} \equiv \sup \left\{ \int_D f^2(x)s(x)\eta(dx) : f \in \mathcal{S}_{m'}, \|f\| = 1 \right\}.$$

Now, by Cauchy-Schwarz $\int_D f^2(x)s(x)\eta(dx) \leq \|f\|_\infty\|s\|$, when $\|f\| = 1$, and so the constant $M_{m'}$ above is bounded by $\|s\|\sqrt{D_{m'}}$. In that case, we can use (7.10-i) to obtain

$$\sqrt{2kM_{m'}x'_{m'}} \leq x_2\sqrt{D_{m'}} + \frac{k\|s\|}{2x_2}x'_{m'},$$

for any $x_2 > 0$. On the other hand, by hypothesis $D_{m'} \leq T$, and (7.19) implies that

$$\sqrt{T}\chi_{m'} \leq (1 + x_1)\sqrt{V_{m'}} + x_2\sqrt{D_{m'}} + \left( \frac{k\|s\|}{2x_2} + k(x_1) \right) x'_{m'}.$$

Choosing the constant $x'_{m'}$ as

$$x'_{m'} = \frac{x_3\sqrt{d_{m'}}}{\frac{k\|s\|}{2x_2} + k(x_1)} + \xi,$$

we get that for any $x_1 > 0$, $x_2 > 0$, $x_3 > 0$, and $\xi > 0$,

(7.21) $$\sqrt{T}\chi_{m'} \leq (1 + x_1)\sqrt{V_{m'}} + x_2\sqrt{D_{m'}} + x_3\sqrt{d_{m'}} + f_1(\xi),$$

with probability larger than $1 - \mathcal{K}_1 e^{-\xi}$, where

(7.22) $$f_1(\xi) = \left( \frac{k\|s\|}{2x_2} + k(x_1) \right) \xi,$$

$$\mathcal{K}_1 = \Gamma \sum_{n=1}^\infty n^R \exp \left( -\sqrt{n}x_3 \Big/ \left( \frac{k\|s\|}{2x_2} + k(x_1) \right) \right).$$

Squaring (7.21) and using (7.10-ii) repeatedly, we conclude that, for any $y > 1$, $x_2 > 0$, and $x_3 > 0$, there exist both a constant $\mathcal{K}_1 > 0$ and a quadratic function of the form $f_2(\xi) = a\xi^2$ (independent of $T$, $m'$, and of the family of linear models) such that, with probability greater than $1 - \mathcal{K}_1 e^{-\xi}$,

(7.23) $$\chi^2_{m'} \leq y\frac{V_{m'}}{T} + x_2\frac{D_{m'}}{T} + x_3\frac{d_{m'}}{T} + \frac{f_2(\xi)}{T}, \quad \forall m' \in \mathcal{M}.$$

Then, (7.18) immediately follows from (7.23) and (7.16). $\qquad\square$

*Proof of (4.5) for the case (c).* By the inequality (7.4), we can upper bound $V_{m'}$ by $\hat{V}_{m'}$ on an event of large probability. Namely, for every $x'_{m'} > 0$ and $x > 0$, with probability greater than $1 - \sum_{m' \in \mathcal{M}} e^{-x'_{m'}}$,

$$(7.24) \qquad (1 + x)\left(\hat{V}_{m'} + \left(\frac{1}{2x} + \frac{5}{6}\right)\frac{D_{m'}}{T}x'_{m'}\right) \geq V_{m'}, \;\; \forall m' \in \mathcal{M},$$

(recall that $D_m = \|\sum_{i=1}^{d_m} \varphi_{i,m}^2\|_\infty$). Since by hypothesis $D_{m'} < T$, and choosing

$$x'_{m'} = x'd_{m'} + \xi, \;\; (x' > 0),$$

it is seen that for any $x > 0$ and $x_4 > 0$, there exist a positive constant $\mathcal{K}_2$ and a function $f(\xi) = b\xi$ (independent of $T$ and of the linear models) such that with probability greater than $1 - \mathcal{K}_2 e^{-\xi}$

$$(7.25) \qquad (1 + x)\hat{V}_{m'} + x_4 d_{m'} + f(\xi) \geq V_{m'}, \;\; \forall m' \in \mathcal{M}.$$

Here, we get $\mathcal{K}_2$ from the polynomial assumption on the class of models. Combining (7.25) and (7.18), it is clear that for any $y_2 > 1$, and positive $x_i$, $i = 1, 2, 3$, we can choose a pair of positive constants $C_1 < 1$, $C'_1 > 1$, an increasing quadratic polynomial of the form $f(\xi) = a\xi^2 + b\xi$, and a constant $\mathcal{K} > 0$ (all independent of the family of linear models and of $T$) so that, with probability greater than $1 - \mathcal{K}e^{-\xi}$

$$(7.26) \qquad \begin{aligned} C_1\|s - \tilde{s}\|^2 &\leq C'_1\|s - s_m^\perp\|^2 \\ &+ y_2\frac{\hat{V}_{\hat{m}}}{T} + x_1\frac{D_{\hat{m}}}{T} + x_2\frac{d_{\hat{m}}}{T} - \text{pen}(\hat{m}) \\ &+ x_3\frac{D_m}{T} + \text{pen}(m) + \frac{f(\xi)}{T}. \end{aligned}$$

Next, we take $y_2 = c$, $x_1 = c'$, and $x_2 = c''$ to cancel $-pen(\hat{m})$ in (7.26). By Lemma 7.4, it follows that

$$(7.27) \qquad C_1\mathbb{E}\left[\|s - \tilde{s}\|^2\right] \leq C'_1\|s - s_m^\perp\|^2 + \left(1 + \frac{x_3}{c'}\right)\mathbb{E}\left[\text{pen}(m)\right] + \frac{C''_1}{T}.$$

Since $m$ is arbitrary, we obtain the case (c) of (4.5). □

*Proof of (4.5) for the case (a).* One can bound $V_{m'}$, as given in (7.20), by $D_{m'}\rho$ (assuming that $\rho < \infty$). On the other hand, (7.4) implies that

$$(7.28) \qquad (1 + x_1)\left(\frac{\mathcal{N}}{T} + \left(\frac{1}{2x_1} + \frac{5}{6}\right)\frac{\xi}{T}\right) \geq \rho,$$

with probability greater than $1 - e^{-\xi}$. Using these bounds for $V_{m'}$ and the assumption that $D_{m'} \leq T$, (7.18) reduces to

$$(7.29) \qquad \begin{aligned} C_1\|s - \tilde{s}\|^2 &\leq C'_1\|s - s_m^\perp\|^2 \\ &+ y\frac{D_{\hat{m}}\mathcal{N}}{T^2} + x_1\frac{d_{\hat{m}}}{T} - \text{pen}(\hat{m}) \\ &+ x_2\frac{D_m\mathcal{N}}{T^2} + \text{pen}(m) + \frac{f(\xi)}{T}, \end{aligned}$$

which is valid with probability $1 - \mathcal{K}e^{-\xi}$. In (7.29), $y > 1$, $x_1 > 0$ and $x_2 > 0$ are arbitrary, while $C_1$, $C'_1$, the increasing quadratic polynomial of the form $f(\xi) =$ $a\xi^2 + b\xi$, and a constant $\mathcal{K} > 0$ are determined by $y$, $x_1$, and $x_2$ independently of the family of linear models and of $T$. We point out that we divided and multiplied by $\rho$ the terms $D_{\hat{m}}/T$ and $D_m/T$ in (7.18), and then applied (7.28) to get (7.29). It is now clear that $y = c$, and $x_1 = c'$ will produce the desired cancelation. □

*Proof of (4.5) for the case (b).* We first upper bound $D_{\hat{m}}$ by $\beta^{-1}V_{\hat{m}}$ and $d_{\hat{m}}$ by $(\beta\phi)^{-1}V_{\hat{m}}$ in the inequality (7.18):

$$
\begin{aligned}
(7.30) \qquad C_1\|s - \tilde{s}\|^2 \leq & \; C_1'\|s - s_m^{\perp}\|^2 + \left(y + x_1\beta^{-1} + x_2(\beta\phi)^{-1}\right)\frac{V_{\hat{m}}}{T} \\
& - \mathrm{pen}(\hat{m}) + x_3\beta^{-1}\frac{V_m}{T} + \mathrm{pen}(m) + \frac{f(\xi)}{T}.
\end{aligned}
$$

Then, using $d_{m'} \leq (\beta\phi)^{-1}V_{m'}$ in (7.25) and letting $x_4(\beta\phi)^{-1}$ vary between 0 and 1, we verify that for any $x' > 0$, a positive constant $\mathcal{K}_4$ and a polynomial $f$ can be found so that with probability greater than $1 - \mathcal{K}_4 e^{-\xi}$,

$$
(7.31) \qquad (1 + x')\hat{V}_{m'} + f(\xi) \geq V_{m'}, \quad \forall m' \in \mathcal{M}.
$$

Putting together (7.31) and (7.30), it is clear that for any $y > 1$ and $x_1 > 0$, we can find a pair of positive constants $C_1 < 1, C_1' > 1$, an increasing quadratic polynomial of the form $f(\xi) = a\xi^2 + b\xi$, and a constant $\mathcal{K} > 0$ (all independent of the family of linear models and of $T$) so that, with probability greater than $1 - \mathcal{K}e^{-\xi}$,

$$
\begin{aligned}
(7.32) \qquad C_1\|s - \tilde{s}\|^2 \leq & \; C_1'\|s - s_m^{\perp}\|^2 + y\frac{\hat{V}_m}{T} - \mathrm{pen}(\hat{m}) \\
& + x_1\frac{V_m}{T} + \mathrm{pen}(m) + \frac{f(\xi)}{T}.
\end{aligned}
$$

In particular, by taking $y = c$, the term $-\mathrm{pen}(\hat{m})$ cancels out. Lemma 7.4 implies that

$$
(7.33) \qquad C_1\mathbb{E}\left[\|s - \tilde{s}\|^2\right] \leq C_1'\|s - s_m^{\perp}\|^2 + (1 + x_1)\,\mathbb{E}\left[\mathrm{pen}(m)\right] + \frac{C_1''}{T}.
$$

Finally, (4.5) (b) follows since $m$ is arbitrary. ∎

**Remark 7.5.** Let us analyze more carefully the values that the constants $C$ and $C'$ can take in the inequality (4.5). For instance, consider the penalty function of part (c). As we saw in (7.27), the constants $C$ and $C'$ are determined by $C_1, C_1', C_1''$, and $x_3$. The constant $C_1$ was proved to be $y_1 - 1$ if $1 < y_1 < 2$, while it can be made arbitrarily close to one otherwise (see the comment immediately after (7.16)). On the other hand, $y_1$ itself can be made arbitrarily close to the penalization parameter $c$ since $c = y_2 = y_1(1 + x)y$, where $x$ is as in (7.24) and $y$ is in (7.23). Then, when $c \geq 2$, $C_1$ can be made arbitrarily close to one at the cost of increasing $C_1''$ in (7.27). Similarly, paying a similar cost, we are able to select $C_1'$ as close to one as we wish and $x_3$ arbitrarily small. Therefore, it is possible to find for any $\varepsilon > 0$, a constant $C'(\varepsilon)$ (increasing in $\varepsilon$) so that

$$
(7.34) \qquad \mathbb{E}\|s - \tilde{s}\|^2 \leq (1 + \varepsilon)\inf_{m \in \mathcal{M}}\left\{\|s - s_m^{\perp}\|^2 + \mathbb{E}\left[\mathrm{pen}(m)\right]\right\} + \frac{C'(\varepsilon)}{T}.
$$

A more thorough inspection shows that

$$
\lim_{\varepsilon \to 0} C'(\varepsilon)\varepsilon = K,
$$

where $K$ depends only $c$, $c'$, $c''$, $\Gamma$, $R$, $\|s\|$, and $\|s\|_{\infty}$. The same reasoning applies to the other two types of penalty functions when $c \geq 2$. In particular, we point out that $\widetilde{C}$ can be made arbitrarily close to 2 in the oracle inequality (4.8) at the price of having a large constant $\widetilde{C}'$.

### 7.2. Some additional proofs

*Proof of Corollary 5.1.* The idea is to estimate the bias and the penalized term in (4.5). Clearly, the dimension $d_m$ of $\mathcal{S}_m^k$ is $m(k+1)$. Also, $D_m$ is bounded by $(k+1)^2 m/(b-a)$ (see (7) in [8]), and

$$\mathbb{E}\left[\hat{V}_m\right] = \int_a^b \left(\sum_i \varphi_{i,m}^2(x)\right) s(x) dx \leq (k+1) m \|s\|_\infty,$$

since the functions $\varphi_{i,m}$ are orthonormal. On the other hand, by (10.1) in Chapter 2 of [15], if $s \in \mathcal{B}_\infty^\alpha \left(L^p([a,b])\right)$, there is a polynomial $q \in \mathcal{S}_m^k$ such that

$$\|s - q\|_{L^p} \leq c_{[\alpha]} |s|_{\mathcal{B}_\infty^\alpha(L^p)} (b-a)^\alpha m^{-\alpha}.$$

Thus,

$$\|s - s_m^\perp\| \leq c_{[\alpha]} (b-a)^{\frac{1}{2} - \frac{1}{p} + \alpha} |s|_{\mathcal{B}_\infty^\alpha(L^p)} m^{-\alpha}.$$

By (4.5)), there is a constant $M$ (depending on $C$, $c$, $c'$, $c''$, $\alpha$, $k$, $b-a$, $p$, $|s|_{\mathcal{B}_\infty^\alpha(L^p)}$, and $\|s\|_\infty$), for which

$$\mathbb{E}\left[\|s - \tilde{s}_\tau\|^2\right] \leq M \inf_{m \in \mathcal{M}_T} \left\{m^{-2\alpha} + \frac{m}{T}\right\} + \frac{C'}{T}.$$

It is not hard to see that, for large enough $T$, the infimum on the above right hand side is $O_\alpha(T^{-2\alpha/(2\alpha+1)})$ (where $O_\alpha$ means that the ratio of the terms is bounded by a constant depending only on $\alpha$). Since $M$ is monotone in $|s|_{\mathcal{B}_\infty^\alpha(L^p)}$ and $\|s\|_\infty$, (5.2) is verified. □

*Proof of Lemma 7.1.* Let

$$(7.35) \qquad \gamma_D(f) \equiv -\frac{2}{T} \iint_{[0,T] \times D} f(x)\, \mathcal{J}(dt,dx) + \int_D f^2(x)\, \eta(dx),$$

which is well defined for any function $f \in L^2((D,\eta))$, where $D \in \mathcal{B}(\mathbb{R}_0)$ and $\eta$ is as in (2.6)-(2.8). The projection estimator is the unique minimizer of the contrast function $\gamma_D$ over $S$. Indeed, plugging $f = \sum_{i=1}^d \beta_i \varphi_i$ in (7.35) gives $\gamma_D(f) = \sum_{i=1}^d (-2\beta_i \hat{\beta}_i + \beta_i^2)$, and thus, $\gamma_D(f) \geq -\sum_{i=1}^d \hat{\beta}_i^2$, for all $f \in S$. Clearly,

$$\gamma_D(f) = \|f\|^2 - 2\langle f, s \rangle - 2\nu_D(f) = \|f - s\|^2 - \|s\|^2 - 2\nu_D(f).$$

By the very definition of $\tilde{s}$, as the penalized projection estimator,

$$\gamma_D(\tilde{s}) + \text{pen}(\hat{m}) \leq \gamma_D(\hat{s}_m) + \text{pen}(m) \leq \gamma_D(s_m^\perp) + \text{pen}(m),$$

for any $m \in \mathcal{M}$. Using the above results,

$$\begin{aligned}
\|\tilde{s} - s\|^2 &= \gamma_D(\tilde{s}) + \|s\|^2 + 2\nu_D(\tilde{s}) \\
&\leq \gamma(s_m^\perp) + \|s\|^2 + 2\nu_D(\tilde{s}) + \text{pen}(m) - \text{pen}(\hat{m}) \\
&= \|s_m^\perp - s\|^2 + 2\nu_D(\tilde{s} - s_m^\perp) + \text{pen}(m) - \text{pen}(\hat{m}).
\end{aligned}$$

Finally, notice that $\nu_D(\tilde{s} - s_m^\perp) = \nu_D(\tilde{s} - s_{\hat{m}}^\perp) + \nu_D(s_{\hat{m}}^\perp - s_m^\perp)$ and that $\nu_D(\hat{s}_m - s_m^\perp) = \chi_m^2$. □

*Proof of inequality (7.4).* Just note that for any $a, b, \varepsilon > 0$:

$$(7.36) \qquad a - \sqrt{2ab} - \frac{1}{3}b \geq \frac{a}{1+\varepsilon} - \left(\frac{1}{2\varepsilon} + \frac{5}{6}\right)b.$$

Evaluating the integral in (7.3) for $-f$, we can write

$$\mathbb{P}\left[\int_{\mathrm{X}} f(x)N(dx) \geq \int_{\mathrm{X}} f(x)\mu(dx) - \|f\|_\mu \sqrt{2u} - \frac{1}{3}\|f\|_\infty u\right] \geq 1 - e^{-u}.$$

Using $\|f\|_\mu^2 \leq \|f\|_\infty \int_{\mathrm{X}} |f(x)|\mu(dx)$ and (7.36), lead to

$$\mathbb{P}\left[\int_{\mathrm{X}} f(x)N(dx) \geq \frac{1}{1+\varepsilon}\int_{\mathrm{X}} f(x)\mu(dx) - \left(\frac{1}{2\varepsilon} + \frac{5}{6}\right)\|f\|_\infty u\right] \geq 1 - e^{-u},$$

which is precisely the inequality (7.4). □

*Proof of Lemma 7.4.* Let $Z^+$ be the positive part of $Z$. First,

$$\mathbb{E}[Z] \leq \mathbb{E}[Z^+] = \int_0^\infty \mathbb{P}[Z > x]dx.$$

Since $h$ is continuous and strictly increasing, $\mathbb{P}[Z > x] \leq K\exp(-h^{-1}(x))$, where $h^{-1}$ is the inverse of $h$. Then, changing variables to $u = h^{-1}(x)$,

$$\int_0^\infty \mathbb{P}[Z > x]dx \leq K\int_0^\infty e^{-h^{-1}(x)}dx = K\int_0^\infty e^{-u}h'(u)du.$$

Finally, an integration by parts yields $\int_0^\infty e^{-u}h'(u)du = \int_0^\infty h(u)e^{-u}du.$ □

**Acknowledgments.** The authors are grateful to an anonymous referee for helpful comments and suggestions. It is also a pleasure to thank P. Reynaud-Bouret for helpful discussions.